\documentclass{article}
\usepackage{amsmath,amsthm,amssymb}
\usepackage{graphicx}
\usepackage{psfrag}
\usepackage{booktabs}
\usepackage{subfig}
\usepackage{algorithm}
\usepackage{algorithmic}
\theoremstyle{plain}

\title{
Interpolation methods to compute statistics of a stochastic partial differential equation}

\author{Daniela Steffes-lai (1), Eveline Rosseel (2), Tanja Clees (1) \\
\small
((1) Fraunhofer SCAI, Sankt Augustin, Germany, \\
\small
 (2) Katholieke Universiteit Leuven, Leuven, Belgium )}

\begin{document}
\thispagestyle{empty}
\pagenumbering{arabic}
 \maketitle

\begin{abstract}
This paper considers the analysis of partial differential equations (PDE) containing multiple random variables. Recently developed collocation methods enable the construction of high-order stochastic solutions by converting a stochastic PDE into a system of deterministic PDEs. This interpolation method requires that the probability distribution of all random input variables is known a priori, which is often not the case in industrially relevant applications. Additionally, this method suffers from a curse of dimensionality, i.e., the number of deterministic PDEs to be solved grows exponentially with respect to the number of random variables. This paper presents an alternative interpolation method, based on a radial basis function (RBF) metamodel, to compute statistics of the stochastic PDE. The RBF metamodel can be constructed even if the probability distribution of all random variables is not known. Then, a lot of statistic scenarios with different probability distributions of the random variables can be computed with this single metamodel. In order to reduce the model complexity, we present a parameter screening technique which can be combined with an interpolation method to solve a reduced stochastic model. Numerical results of a model problem demonstrate that the RBF metamodel is as fast as a low order collocation approach and achieves a good accuracy. The parameter screening is able to reduce the dimension and, thus, to accelerate the computation of the stochastic solution.
\end{abstract}


\section{Introduction}
During the fabrication of products important material and process parameters, geometries and also external influences can vary considerably. These variations can have a substantial influence on the quality of the resulting products. The current trend goes to more accuracy and therefore more robust predictions. 
Hence, there is a high need for efficient stochastic computations, combined with an efficient analysis of sensitivity and robustness, see, e.g., the review papers~\cite{Xiu2009, Matthies2008, Najm2009, Sudret2008}.
The objective is to model uncertainty from the beginning of the simulation, and not as criterion afterwards.
Mathematically, parameters which are known only approximately can be modeled as stochastic variables or processes. The model takes the form of a partial differential equation (PDE) with coefficients that are random variables or random fields, e.g., in the differential operator, the initial or boundary conditions. These random variables are mostly not well known in practice.

We seek for solutions of a stochastic PDE, called stochastic solution, that contain all possible probability information, including a complete description of the probability density function (pdf). Until recently Monte Carlo simulations were the standard tool for solving stochastic PDEs. Today, new stochastic finite element methods~\cite{Xiu2009,Matthies2008} are becoming popular because they enable the computation of high-order stochastic solutions with substantial less computational effort than Monte Carlo simulations. A prominent interpolation variant is the collocation method~\cite{Babuska2007,Mathelin2005, Xiu2005}.  
This approach transforms the stochastic PDE into a set of decoupled deterministic PDEs. However, this method suffers from the ``\textit{curse of dimensionality}''~\cite{bellman61, Donoho2000}, that is, the computational cost grows exponentially as a function of the number of random variables.

The main focus of this work is the computation of statistics of a stochastic PDE by means of interpolation methods. We investigate an alternative approach to the collocation method, namely a radial basis function (RBF) metamodel, which is accelerated by a singular value decomposition (SVD).
The difference to the collocation approach are mainly the construction of the sample points and the interpolation function itself. The accelerated metamodel does not require the probability distribution of the random input variables in the setup phase. The probability distribution of the random input variables is only used to compute statistics of the solution.

Additionally, we perform a sensitivity analysis, which determines the dependence of the solution on the variation of random variables, also called uncertain \textit{design-parameters}. One goal of a sensitivity analysis is to identify the design-parameters possessing the most influence on the solution. This subject is often also called \textit{parameter screening}. We show that the parameter screening approach and the accelerated metamodel approach can be combined in order to reduce the curse of dimensionality.

We investigate the computation of statistics of the solution. We are interested in quantifying the \textit{robustness} of the solution, that is, a small change in the design-parameters should also result in only a small change of the solution. Therefore, we will focus on the median and the $0.68$-quantile of the solution. We are not primarily interested in the tails of the pdf, since these represent the mostly improbable cases, or worst case scenarios. The tails of the pdf are investigated in a reliability analysis, which is beyond our scope.

This paper is organized as follows. Section \ref{sec:sPDE} introduces the stochastic partial differential equation to be solved. Section \ref{sec:sensan} presents a simple and effective approach for parameter screening, together with a discussion of its advantages and limitations.
The collocation approach and a RBF metamodel to compute statistics of the solution are presented in Section~\ref{sec:Interpolation}. Section~\ref{sec:Results} describes the model problem for numerical tests. Numerical results of the RBF metamodel approach are given and compared to those of the collocation method. We conclude the paper with a discussion and some future work issues in Section~\ref{sec:concl}.

\section{Problem definition}\label{sec:sPDE}
Consider a stochastic partial differential equation (sPDE)
\begin{align}\label{eq:model}
 \mathcal{L}(\boldsymbol{x}, \omega; u(\boldsymbol{x},\omega)) &= b(\boldsymbol{x},\omega) \qquad \boldsymbol{x} \in \mathbf{D},\: \omega\in \Omega, \\ \nonumber
u(\boldsymbol{x},\omega) &= g(\boldsymbol{x},\omega) \qquad \boldsymbol{x} \in \partial\mathbf{D}_D,\: \omega\in \Omega, \\
\frac{\partial u(\boldsymbol{x},\omega)}{\partial n} &= h(\boldsymbol{x},\omega) \qquad \boldsymbol{x} \in \partial\mathbf{D}_N,\: \omega\in \Omega, \nonumber
\end{align}
where $\mathcal{L}$ represents a general -- possibly nonlinear -- differential operator containing stochastic coefficients, $\mathbf{D}$ is the spatial domain enclosed by a Dirichlet boundary $\partial\mathbf{D}_D$ and a Neumann boundary $\partial \mathbf{D}_N$, and $\Omega$ is a sample space. 

We make use of the finite dimensional noise assumption~\cite{Babuska2007}. This assumption states that all sources of randomness can be approximated by a finite (small) number of random variables. This is, e.g., satisfied when the random processes present in the problem have the form of a truncated Karhunen-Lo\`eve expansion~\cite{Loeve1977}.
We denote the vector of all random variables in \eqref{eq:model} by $\boldsymbol{\xi}:= (\xi_1(\omega), \hdots,\xi_L(\omega))$, with $L$ the total number of random variables. Each random variable~$\xi_i(\omega):\Omega \rightarrow \Gamma_i $ is a function from the event space $\Omega$ onto a part of the real axis, $\Gamma_i \subset \mathbb{R}$, and is characterized by its probability density function $\rho_i(y_i)$, with $y_i \in \Gamma_i$. Assuming that all random variables are independent, their joint probability distribution equals $\rho(\boldsymbol{y}) = \prod_{i=1}^L \rho_i(y_i)$, $\boldsymbol{y}=(y_1,\hdots,y_L)\in \Gamma =\prod_{i=1}^L \Gamma_i$. Using these assumptions, we can rewrite the sPDE~\eqref{eq:model} as a parametric deterministic PDE
\begin{align}\label{eq:modelRandom}
  \mathcal{L}(\boldsymbol{x}, \boldsymbol{y}; u(\boldsymbol{x},\boldsymbol{y})) &= b(\boldsymbol{x},\boldsymbol{y}) \qquad \boldsymbol{x} \in \mathbf{D},\: \boldsymbol{y} \in \Gamma, \\ \nonumber
u(\boldsymbol{x},\boldsymbol{y}) &= g(\boldsymbol{x},\boldsymbol{y}) \qquad \boldsymbol{x} \in \partial\mathbf{D}_D,\: \boldsymbol{y} \in \Gamma, \\
\frac{\partial u(\boldsymbol{x},\boldsymbol{y})}{\partial n} &= h(\boldsymbol{x},\boldsymbol{y}) \qquad \boldsymbol{x} \in \partial\mathbf{D}_N,\: \boldsymbol{y} \in \Gamma. \nonumber
\end{align}
In the context of the parameter screening technique presented in the next section, we will call the random variables $\xi$ \textit{design-parameters}. The space spanned by these design-parameters is called \textit{design-parameter space}.

Note, if \eqref{eq:modelRandom} is solved on a spatial discretization with $M$ degrees of freedom, the discretized solution in node $i$ is denoted by
\begin{equation}\label{eq:solutionDiscrete}
u_i(\boldsymbol{y}), \quad(i=1,\dots,M).
\end{equation}

\section{Parameter screening}\label{sec:sensan}
We present a \textit{parameter screening} method that allows the identification of the design-parameters possessing the most influence on the solution. This is used to reduce the dimension of the design-parameter space and therefore the model complexity.
To investigate the model and the influences of variations of design-parameters onto the solution, a special design of experiments is performed and described in Section~\ref{sec:sensan:Doe}. Section~\ref{sec:sensan:JH} describes the sensitivity analysis in detail. We start with the computation of an indicator for nonlinearity. Then we compute a measure which also deals with the linear behavior of design-parameters. The advantages and limitations of the presented parameter screening method are discussed in Section~\ref{sec:sensan:limit}.
\subsection{Design of experiments}\label{sec:sensan:Doe}
First, a small design of experiments (DoE)~\cite{Myers2009} is generated, consisting only of $2L+1$ sampling points $\{\boldsymbol{\zeta}_1,\dots,\boldsymbol{\zeta}_{2L+1}\}$, with $L$ the total number of random variables. The sampling points represent variations of the random variables $\xi$, often so called \textit{design-parameters}. Each random variable is varied separately. The first sample corresponds to the mean values of all random variables. In addition, $2L$ samples are generated by sampling each random variable in a minimum and maximum value, while the other random variables are sampled in their mean value. 

Next, deterministic simulations on these sampling points are carried out to compute discrete solutions $u_i(\boldsymbol{\zeta}_1),\hdots, u_i(\boldsymbol{\zeta}_{2L+1}), i = 1,\dots, M $ on a spatial discretization with $M$ spatial degrees of freedom. We define a discrete, deterministic solution vector as $\boldsymbol{u}(\boldsymbol{\zeta}_j) = (u_1(\boldsymbol{\zeta}_j),\dots,u_M(\boldsymbol{\zeta}_j))^T$. All solution vectors are grouped together column by column in a matrix $\mathbf{X} \in \mathbb{R}^{M \times (2L+1)}$. 
\begin{equation}\label{eq:database}
\mathbf{X} := \mathbf{X(\boldsymbol{\zeta})}= \left[ \boldsymbol{u}(\boldsymbol{\zeta}_1)\hdots \boldsymbol{u}(\boldsymbol{\zeta}_{2L+1}) \right].
\end{equation}
The rows of $\mathbf{X}$ correspond to the $M$ spatial unknowns. $\mathbf{X}$ is called \textit{data base} for the following analysis. Each column of the data base represents a solution to one particular parameter variation.
\subsection{Sensitivity analysis}\label{sec:sensan:JH}
A sensitivity analysis determines the dependence of the solution on the variation of design-parameters. It allows one to identify the design-parameters possessing the most influence on the solution \cite{Hamby1994, Rabitz1983, Saltelli2008}. This objective of sensitivity analysis is often called \textit{parameter screening}. Design-parameters possessing no or only a small influence on the solution are negligible and can be fixed to their mean value. This leads to a reduction of the dimension of the design-parameter space and therefore to a reduction of the model complexity. Sensitivity methods can be divided into local and global methods, for a detailed discussion refer to \cite{Rabitz1983, Saltelli2008}. A simple local sensitivity method is to approximate partial derivatives of the solution w.r.t. the design-parameters with a finite difference scheme. We apply this method on the data base \eqref{eq:database}. 

We compute $\boldsymbol{J}=(J_{ij})_{i=1,\dots,M, j=1,\dots,L}$ with a finite difference scheme as a second order approximation of the Jacobian matrix, where 
 $J_{ij}$ is the sensitivity measure for design-parameter $j$ in grid point $i$.
We use
\begin{equation}\label{eq:sensitivityLinGlob}
S_j = \sqrt{\sum_{i=1}^M J_{ij}^2}
\end{equation}
as global measure on the whole grid for estimating the linear influence of the $j$th parameter. 
The design-parameters are ranked by sorting the sensitivity measures $S_j, j=1,\dots,L$ in descending order. The dimension of the design-parameter space is then reduced by retaining only the most sensitive parameters in the model.

The Jacobian matrix only shows the magnitude of the linear influence of each design-parameter onto the solution. To overcome this drawback, additionally an indicator for the nonlinear influence of the design-parameters on the solution is computed. This indicator makes use of the diagonal part of the Hessian matrix; the full Hessian matrix is not computed in order to keep the computation time to a minimum. The diagonal part of the Hessian matrix $\boldsymbol{\mathrm{diag}H}=(\mathrm{diag}H_{ij})$ is analogously approximately computed with a finite difference scheme of second order. 
We use
\begin{equation}\label{eq:sensitivityNonlinGlob}
S2_j = \sqrt{\sum_{i=1}^M \mathrm{diag}H_{ij}^2}
\end{equation}
as global measure of the nonlinear influence of the $j$th parameter.
We assume that all parameters have the same input range. This can easily be achieved by a normalization of the parameter input space, e.g. to $[-1,1]$. Then, we are allowed to compare the measures for estimating the linear as well as nonlinear influence of different parameters.
If
\begin{equation}\label{eq:criteriumGlob}
S_j < c S2_j,
\end{equation} 
with $c$ appropriately chosen, e.g. $c = \sigma_j$, the dependence of the $j$th design-parameter on the solution should be assumed nonlinear. Additionally, local effects can be analysed by pointwise comparing the entries of the Jacobian and the diagonal part of the Hessian matrix, i.e., if
\begin{equation}\label{eq:criteriumLoc}
|J_{ij}| <c |\mathrm{diag}H_{ij}|,
\end{equation}
the dependence of $j$th design-parameter on the solution in grid point $i$ should be assumed to be nonlinear.
The opposite is not true. It is possible that the Hessian matrix is not diagonal dominant, especially it can occur that $\boldsymbol{\mathrm{diag}H} = 0$, but the off-diagonal elements are large.

In order to ensure that the design-parameters behave linear on the solution (i.e. at least the derivatives of second order are small), the full Hessian matrix has to be analysed. This requires the computation of $L(L-1)/2$ additional entries of the Hessian matrix due to its symmetry. Thus, we need to perform a DoE with $4 L(L-1)/2$ additional simulations to approximate the required second order derivatives w.r.t. the $j$th and $l$th parameter in grid point $i$ with finite difference schemes of second order. 
Analogously to the diagonal part, we define a global approximation of the $j,l$th entry of the Hessian matrix. We define a global measure $D$ \cite{Nikitin2012} as
\begin{equation}\label{eq:measureD}
D:= \sqrt{\sum_{j=1}^L J_j} - \frac{\sigma}{4} |\alpha_{\mathrm{max}}|,
\end{equation}
where $\alpha_{\mathrm{max}}$ denotes the maximal eigenvalue of the global approximation of the Hessian matrix and $\sigma$ is the maximum of $\sigma_j, j=1, \dots, L$.
If $D$ is negative, the dependence of design-parameters on the solution is assumed to be nonlinear, otherwise it is linear.

In summary, the sensitivity analysis computes the linear sensitivity measure $\boldsymbol{J}$ and the indicator for nonlinearity $\boldsymbol{\mathrm{diag}H}$. Only, if the inequalities \eqref{eq:criteriumGlob} and \eqref{eq:criteriumLoc} do not hold, the full Hessian matrix is computed.
A reduced stochastic model will be constructed by replacing nearly-linear design-parameters which have only a small influence on the solution, by their mean value.
\subsection{Discussion}\label{sec:sensan:limit}
\subsubsection{Advantages} 
Firstly, the presented DoE-based parameter screening method is independent of the probability distribution of the random variables. This makes the approach very general and useful in industrial relevant applications, because usually the correct probability distribution of the design-parameters are not known.

Secondly, the DoE needs only $2L+1$ samples, this is the minimum number of samples needed to compute a second order approximation of the Jacobian as well as the diagonal part of the Hessian matrix. This number of samples is appropriate for representing a linear dependence of each parameter on the solution. As a result, the local sensitivity analysis by means of finite differences requires less computational effort than variance based sensitivity methods, e.g., Sobol indices, for which many simulations are to be performed in Monte Carlo alike algorithms \cite{Saltelli2008, Sobol2001}. Additionally, the indicator of nonlinear influence gives an important hint to the nature of design-parameter behavior. If the number of simulations, which can be performed, is not restricted by the simulation time, additionally, the full Hessian matrix and the measure $D$ is computed. This approach ensures that only parameters with a small (the last few parameters in the ranking), linear influence are neglected.

Finally, the presented approach for sensitivity analysis is able to locate local areas on the grid, where the behavior of design-parameters onto the solution is nonlinear. Hence, when solving a sPDE with an interpolation based method, more advanced interpolation methods could be applied in the regions of interest in order to improve the accuracy of the interpolation result.

\subsubsection{Limitations}
Sensitivity measures computed with partial derivatives around a nominal value are only valid if a very low order Taylor series approximation of the solution is valid. This is typically true only in small regions around the nominal value, i.e., if $\sigma_j$ is chosen appropriately small. This drawback is somewhat captured with the computation of the indicator of nonlinearity. Note, that interactions of third and higher order can only be computed with more computational effort.
Secondly, interactions among design-parameters can not be detected  by varying only one parameter at a time. Therefore, this approach requires independent random variables. This limitation can be avoided by computing the full Hessian matrix.

\section{Interpolation methods}\label{sec:Interpolation}
Let an unknown or  computational expensive function $u$ be determined by a set of $N$ sample points $\{\zeta_i\}_{i=1,\dots,N}$ in the parameter space and corresponding solutions $\{u(\zeta_i)\}_{i=1,\dots,N}$.
We seek an approximant $\tilde{u}$ for $u$, which is faster to evaluate. Additionally, we require
\begin{equation}\label{eq:interpolation}
\tilde{u}(\zeta_i) = u(\zeta_i),\ i=1,\dots,N,
\end{equation}
i.e., the discrete sampling points are interpolated. Note, we write $u_i$ for $u(\zeta_i)$.
There are a lot of commonly used interpolation methods, for example \textit{metamodels}, also called \textit{response surface models}, which are  widely used in the field of optimization (for an overview see \cite{Myers2009, Fang2005, Kleijnen2005}). Note, for some of these approximation methods, the stronger interpolation condition~\eqref{eq:interpolation} does not hold.
Additionally, there exist special interpolation methods for solving sPDEs, for example the stochastic collocation method. We will use two interpolation methods, the stochastic collocation method and a radial basis function (RBF) metamodel, to compute statistics of the stochastic solution of~\eqref{eq:model}. The stochastic collocation approach is described in Section~\ref{ssec:coll}. Section~\ref{sec:RBFModel} presents metamodels with radial basis functions accelerated by a principal trends analysis.
\subsection{Stochastic collocation finite element method}\label{ssec:coll}
The stochastic collocation finite element method~\cite{Babuska2007,Xiu2005} is a sampling based solution method for stochastic PDEs. In contrast to the Monte Carlo simulation method~\cite{Fishman1996}, it achieves a fast convergence rate by a good choice of the multidimensional stochastic collocation points.
\subsubsection{Collocation approach}
The stochastic collocation method represents the solution discretely by an expansion with multivariate Lagrange polynomials $l_k(\boldsymbol{y})$,
\begin{equation}\label{eq:lagrange}
  u(\boldsymbol{x}, \boldsymbol{y}) \approx \sum_{k=1}^{N_c} u(\boldsymbol{x},\zeta_k) l_k(\boldsymbol{y}).
\end{equation}
The Lagrange polynomials~$l_k(\boldsymbol{y})$ are interpolatory polynomials, defined by a set of multidimensional collocation points,~$\{\zeta_1,\hdots,\zeta_{N_c}\}$. Each sample or collocation point~$\zeta_k$ consists of $L$ components, $\zeta_k = (\zeta_{k,1},\hdots,\zeta_{k,L})$, according to the $L$ random variables present in the problem. 

In order to obtain a fast convergence of~\eqref{eq:lagrange}, the collocation points $\zeta_k$ are typically constructed as a sparse grid of Gauss or Clenshaw-Curtis points.
By using sparse grids, the curse of dimensionality associated with full tensor product collocation grids is reduced. An extended overview of the stochastic collocation method can be found in \cite{Babuska2010, Back2009}.

The stochastic collocation method proceeds by applying a collocation step, i.e., it requires that the residual vanishes at each collocation point~$\zeta_k$,
\begin{equation}\label{eq:collocation}
  \mathcal{L}\left(\boldsymbol{x},\zeta_k; \sum_{m=1}^{N_c} u(\boldsymbol{x},\zeta_m)l_m(\zeta_k) \right) - b(\boldsymbol{x},\zeta_k) = 0 \qquad \forall k = 1,\hdots,N_c.
\end{equation}
Since Lagrange polynomials satisfy the property that $l_k(\zeta_m) = \delta_{k,m}$, with $\delta_{k,m}$ the Kronecker delta, 
\eqref{eq:collocation} corresponds to a set of $N_c$ decoupled deterministic PDEs,
\begin{equation}\label{eq:coll2}
 \mathcal{L}\left(\boldsymbol{x},\zeta_k;  u(\boldsymbol{x},\zeta_k) \right) - b(\boldsymbol{x},\zeta_k) = 0 \qquad \forall k = 1,\hdots,N_c.
\end{equation}
These deterministic PDEs~\eqref{eq:coll2} are then discretized in the spatial domain; we will denote by $M$ the corresponding number of spatial degrees of freedom.

\subsubsection{Computation of statistics}
From the deterministic solutions of \eqref{eq:coll2}, $u(\boldsymbol{x},\zeta_1),\hdots, u(\boldsymbol{x},\zeta_{N_c})$, the statistics of the solution can be obtained. By choosing the collocation points to be a cubature point set~\cite{Xiu2009,Xiu2005}, the mean~$\overline{u}(\boldsymbol{x}) $ and variance $\sigma^2_{u}(\boldsymbol{x})$ can be computed as~\cite{Loeven2007}
\begin{equation*}
\overline{u}(\boldsymbol{x}) = \sum_{k=1}^{N_c} u(\boldsymbol{x},\zeta_k)w_k \qquad \mathrm{and} \qquad \sigma^2_{u}(\boldsymbol{x}) = \sum_{k=1}^{N_c} u(\boldsymbol{x},\zeta_k)^2 w_k - \overline{u}(\boldsymbol{x})^2,
\end{equation*}
where $w_k$, $k= 1,\hdots, N_c$, are the cubature weights corresponding to the cubature points $\zeta_k$.
Other statistics, e.g., quantiles or the probability density function, can be straightforwardly computed after sampling the stochastic solution~\eqref{eq:lagrange} in the $\boldsymbol{y}$-parameters by a Monte Carlo procedure.

\subsubsection{Discussion}
\paragraph{Advantages}
As a sampling-based method, the stochastic collocation approach permits a direct reuse of deterministic simulation software. As a consequence, any nonlinearity present in the stochastic PDE does not increase the complexity of the method.

The stochastic collocation method enables one to compute high-order stochastic solutions with a limited number of deterministic simulations. For stochastic elliptic problems, theoretical convergence estimates in \cite{Babuska2007,Babuska2010,Nobile2008} indicate an algebraic convergence rate of the solution w.r.t.~the number of collocation points, and an exponential rate w.r.t.~the polynomial degree in the case of a tensor product collocation point grid of Gauss points.

To increase the accuracy of a stochastic collocation solution, an extended set of collocation points is required. In the case of nested 1D point sets, e.g., when using Clenshaw-Curtis points, previously obtained deterministic solutions can be reused and only a limited number of additional deterministic simulations need to be performed.
\paragraph{Limitations}
The number of deterministic simulations, i.e., the number of collocation points, increases proportionally to the number of random variables and the required accuracy, which determines the number of 1D collocation points in the interpolation formulae of tensor and sparse grid collocation points. Although sparse grids try to eliminate the curse of dimensionality present in the tensor product approach, the total number of collocation points still grows exceedingly fast in the case of a large number of random variables. Current research aims at further reducing the required number of simulations, for example by taking the anisotropy of the problem into account~\cite{Nobile2008a,Bieri2009}.

\subsection{A RBF metamodel accelerated by a principal trends analysis}\label{sec:RBFModel}
We present a nonlinear interpolatory metamodel with radial basis functions accelerated by a principal trends analysis. First, we introduce RBF metamodels in Section~\ref{sec:RBF} and explain principal trends analysis briefly in Section~\ref{sec:svd}. Then, we describe the algorithm of a RBF metamodel accelerated by a principal trends analysis in Section~\ref{sec:Algorithm}. Finally, in Section~\ref{sec:stats}, we show how the algorithm can be used to compute statistics.
\subsubsection{RBF metamodels}\label{sec:RBF}
A radial basis function $\phi: \mathbb{R} \rightarrow \mathbb{R}$ is a real-valued function whose values depend only on the distance to the origin, $\phi(x) = \phi(||x||)$, where $||\cdotp ||$ denotes the Euclidean norm. Table~\ref{tab:rbf} shows some commonly used types of RBFs for $r=||x||$, \cite{Buhmann2003, Carr2001, Wendland2005}.
\begin{table}[htbp]
\centering
\caption{Commonly used types of radial basis functions}
\renewcommand{\arraystretch}{1.5}
\setlength{\arrayrulewidth}{0.8pt}
\begin{tabular}{|c|c|}
\hline $\phi(r)=r$ & linear \\
\hline $\phi(r)=r^2\log{r}$ & thin-plate spline \\
\hline $\phi(r)=\exp{(-\gamma r^2)}, \, \gamma>0$ & Gaussian \\
\hline $\phi(r)= \sqrt{r^2+c^2}, \, c>0$ & Multiquadric \\
\hline $ \phi(r)=r^3$ & triharmonic spline \\
\hline $ \phi(r)=(r^2+c^2)^{-\frac{1}{2}}$ & inverse Multiquadric \\
\hline 
\end{tabular}
\label{tab:rbf}
\end{table}
A RBF metamodel is a linear combination of radial basis functions \cite{Buhmann2003},
\begin{equation}\label{eq:rbf-metamodel}
\tilde{u}(\zeta) = \sum_{j=1}^N \phi(||\zeta-\zeta_j||) c_j,
\end{equation}
with coefficients $c_i \in \mathbb{R}$, which are determined so that \eqref{eq:interpolation} holds, i.e.,
\begin{equation}
\tilde{u}(\zeta_i) = u_i = \sum_{j=1}^N \phi(||\zeta_i-\zeta_j||) c_j, \quad i=1,\dots,N.
\end{equation}
We define an interpolation matrix $\Phi = (\phi(||\zeta_i-\zeta_j||))_{i,j=1,\dots,N}$. Note that the interpolation matrix is always nonsingular for the multiquadric function, details on invertibility of the interpolation matrix can be found in~\cite{Buhmann2003, Wendland2005}.
By inverting the interpolation matrix, the result $\tilde{u}$ can be written as weighted sum \cite{Nikitina11}
\begin{equation}\label{eq:weightsSum}
\tilde{u}(\zeta) = \sum_{i=1}^N w_i(\zeta)u_i, 
\end{equation}
with 
\begin{equation}\label{eq:weights}
w_i(\zeta) = \sum_{j=1}^N \Phi^{-1}_{ij}\phi(||\zeta-\zeta_j||),
\end{equation}
where $\Phi^{-1}_{ij}$ are the entries of $\Phi^{-1}$.
This metamodel can be extended by polynomial detrending, that is by adding a polynomial part
\begin{equation}\label{eq:detrending}
\tilde{u}(\zeta) = \sum_{j=1}^N \phi(||\zeta-\zeta_j||) c_j + P(\zeta).
\end{equation}
Polynomial detrending improves the precision of interpolation and leads to an exact representation of polynomials up to a certain degree, for more details see, e.g.,~\cite{Buhmann2003}.
\subsubsection{Principal trends analysis}\label{sec:svd}
Principal component analysis (PCA) is a well known mathematical method for data analysis which provides the principal components of the data, i.e., the dominating dependencies of design-parameters on the solution.
PCA is a linear method, that is, it assumes that the data result from a linear transformation of variables. 

One way to compute the principal components is by applying a singular value decomposition (SVD)~\cite{Lee2007}.
This decomposition exists for each $M \times N$ matrix $\mathbf{X}$, $M \geq N$ and is defined as 
\begin{equation}\label{eq:svd}
\mathbf{X}=\mathbf{F}\mathbf{\Lambda} \mathbf{V}^T,
\end{equation}
where $\mathbf{F} = (f_1,f_2,\hdots,f_N) \in \mathbb{R}^{M \times N}$, $\mathbf{V}^T  = (v_1^T,v_2^T,\hdots,v_N^T) \in \mathbb{R}^{N \times N}$ are orthogonal matrices such that $\mathbf{F}^T\mathbf{F} = \mathbf{I}_N$ and $\mathbf{V}^T\mathbf{V}=\mathbf{V}\mathbf{V}^T=\mathbf{I}_N$.   $\mathbf{\Lambda} \in \mathbb{R}^{N \times N}$ is a diagonal matrix with diagonal elements $\lambda_1 \geq \lambda_2 \geq \dots \geq \lambda_N \geq 0.$ 
The number of non-zero singular values determines the rank $r$ of $\mathbf{X}$, that is $\lambda_{r+1}=\lambda_{r+2} = \hdots = \lambda_N=0.$ The entries $\{\lambda_j\}$ are called the singular values of $\mathbf{X}$.

A rank-$k$ approximation $\mathbf{X}_k$ of $\mathbf{X}$ is defined as
\begin{equation}\label{eq:rank-k-approx}
\mathbf{X_k} = \sum_{i=1}^k f_i \lambda_i v_i^T, \qquad \mathrm{with} \, k<r.
\end{equation}
Then it holds that $\mathbf{X}_k$ is the best rank-$k$ approximation in the Frobenius norm $||\cdot ||_\mathcal{F}$ to the matrix $\mathbf{X}$~\cite{Berry1992}, that is
\begin{equation}\label{eq:approx-err}
\min_{\mathrm{rank}(\mathbf{Y})=k} ||\mathbf{X}-\mathbf{Y}||_\mathcal{F}^2 = ||\mathbf{X}-\mathbf{X}_k||_\mathcal{F}^2 = \lambda_{k+1}^2 + \dots + \lambda_r^2.
\end{equation}
Therefore, a linear dimension reduction with the SVD proceeds by setting the singular values $\lambda_{k+1}$ to $\lambda_r$ to zero. The approximation error of the resulting reconstruction equals
\begin{equation}\label{eq:err}
\mathrm{err}^2 = \sum_{l=k+1}^{N} \lambda_l^2,
\end{equation}
which can be controlled over $k$. Therefore, the SVD gives us a method to approximate the intrinsic dimension and to project the data to that lower dimensional space.

The computation of the SVD of a matrix $\mathbf{X} \in \mathbb{R}^{M \times N}$ is generally based on the eigensystem of either $\mathbf{X}^T\mathbf{X}$ or $\mathbf{X}\mathbf{X}^T$~\cite{Press, Golub1996, Berry2003}.
An efficient implementation of the SVD for $M\gg N$ operates block by block so that it is not needed to store the entire matrix in memory. It is based on a computation of the Gram matrix, computing its spectral decomposition, and finally, the $k$ relevant columns of $\mathbf{F}$~\cite{Nikitina11}. We call this approach \textit{fastSVD}.

\subsubsection{RBF metamodel accelerated by principal trends analysis}\label{sec:Algorithm}
Algorithm~\ref{algo:accMetamodel} describes the procedure to accelerate the RBF metamodel~\eqref{eq:detrending} with a principal trends analysis and to use this for the approximation of a $q-$quantile.
\begin{algorithm}
\caption{Approximation of a $q-$quantile with the accelerated metamodel} \label{algo:accMetamodel}
\begin{algorithmic}[1]
\REQUIRE $\boldsymbol{\xi}, \mathrm{threshold}, q$
\STATE $\{\zeta_i\}_{i=1,\dots,L_S} = \mathrm{DoE(\boldsymbol{\xi})}, \quad \mathbf{X} := \mathbf{X(\boldsymbol{\zeta})}= \left[ \boldsymbol{u}(\boldsymbol{\zeta}_1)\hdots \boldsymbol{u}(\boldsymbol{\zeta}_{2L+1}) \right].$\label{algo:metamodel:simulate}
\STATE $\{w_i\}_{i=1;\dots,N} = \mathrm{metamodel}(\boldsymbol{\zeta}), \quad [\mathbf{F},\mathbf{\Lambda},\mathbf{V}^T] = \mathrm{fastSVD}(\mathbf{X},\mathrm{threshold})$\label{algo:metamodel:setup}
\STATE $ \{\tau_i\}_{i=1,\dots,L_{add}} = \mathrm{DoE(\boldsymbol{\zeta})}$\label{algo:metamodel:q-start}
\FOR{$j=1$ \TO $L_{add}$}
\STATE $\boldsymbol{\tilde{u}}(\tau_j) \leftarrow \mathrm{acceleratedMetamodel}(\tau_j,\boldsymbol{w},\mathbf{F},\mathbf{\Lambda},\mathbf{V}^T)$	\label{algo:metamodel:Mixer}
\STATE $\mathrm{q-quantile} \leftarrow \mathrm{P^2_q}(\boldsymbol{\tilde{u}},q)$ \label{algo:metamodel:quantil} 
\ENDFOR
\end{algorithmic}
\end{algorithm}
First, a design of experiments (DoE) is generated to construct a data base $\mathbf{X}$. Next, a metamodel is set up and a singular value decomposition (SVD) of the database $\mathbf{X}$ is computed. Note, that only the weights \eqref{eq:weights} of the metamodel are computed. These steps, line \ref{algo:metamodel:simulate} and line \ref{algo:metamodel:setup}, are called the setup phase to construct the metamodel. In the setup phase the probability distribution of all random input variables is not required.
Then the SVD of $\mathbf{X}$ is used instead of the sample results $\mathbf{X}$ itself for the construction of the accelerated metamodel \cite{Nikitin2010}. That is, the SVD is combined with the metamodel weights, see line \ref{algo:metamodel:Mixer}, in order to approximate the solution $u$ for a given parameter set $\tau_j$. This procedure can be written as
\begin{equation}\label{eq:mixing}
\tilde{u}(\tau_j) \underset{\eqref{eq:weightsSum}}{=} \sum_{i=1}^N u_i w_i(\tau_j) = \sum_{i=1}^k  f_i \lambda_i v_i^T  w_i(\tau_j)
\end{equation}
The so-constructed metamodel is evaluated in $L_{add}$ sampling points in order to compute statistics of the solution by a quasi Monte Carlo approach. The optimal number of sampling points should be large enough and depends on the model problem. In quasi Monte Carlo methods low discrepancy sequences are used as sampling points which lead to a better performance as standard Monte Carlo methods \cite{Niederreiter1978, Morokoff1995}.
The computation of statistics is performed in lines \ref{algo:metamodel:q-start} to \ref{algo:metamodel:quantil} and is described in the next section.
\subsubsection{Computation of statistics}\label{sec:stats}
To compute statistics of $u$, any method can be used that is based only on samples of the design-parameters. The $q$-quantile for a random variable is defined as the value $x$ such that the probability that the random variable will be less than $x$ is at most $q$, with $q \in [0,1]$. For example, the $0.9$-quantile is the value $x$ such that the probability that the random variable will be less than $x$ is at most 90\%. A $q$-quantile is generally estimated as a weighted sum of the order statistics, with a high number of samples, see, e.g.,~\cite{David2003}. The problem of such approaches is that all samples must be stored.
Therefore, we use the $P^2$ algorithm for dynamic calculation of quantiles, instead of the order statistics, to compute the median and $q$-quantiles of the stochastic solution of~\eqref{eq:modelRandom} in Algorithm~\ref{algo:accMetamodel} (line \ref{algo:metamodel:quantil}). The $P^2$ algorithm by Jain and Chlamtac stores only five markers, which are updated as more samples are generated. The five markers are the minimum, the $q/2$-, $q$- and $(q+1)/2$-quantiles and the maximum. The markers are adjusted with a parabolic prediction ($P^2$) formula. Details on the $P^2$ algorithm can be found in \cite{Jain1985}.
The best number of sampling points depends on the model problem; usually a high number of samples ($L_{add} > 1000$) are required.

In summary, to compute statistics of the solution of a sPDE, we perform a DoE to get $L_{add}$ sampling points according to the probability distribution of the random variables, see Algorithm \ref{algo:accMetamodel} line \ref{algo:metamodel:q-start}. For each sampling point we approximate the corresponding solution with the accelerated metamodel, see line \ref{algo:metamodel:Mixer}. Note that no additional deterministic solutions of the sPDE by means of simulation are needed. Finally, we update the $P^2$ algorithm estimator of the $q$-quantile (line \ref{algo:metamodel:quantil}) in each grid node.
\subsubsection{Discussion}
\paragraph{Advantages}
An advantage of the accelerated RBF metamodel method lies in its flexibility. It does not require the probability distribution of the random input variables in the setup phase. 
As for the collocation method, the accelerated metamodel is non-intrusive. Simulation code available for the deterministic problem can be reused directly.
The choice of the RBF ensures the invertibility of the interpolation matrix, and thus the stability of the method. 
Secondly, the acceleration by the singular value decomposition decreases the computational time considerably. The SVD computes the principal trends in the data, and, therefore, takes the correlation between the nodal points into account. Due to this reason, fewer sample points may be necessary in order to get the same accuracy as with using a collocation method. Additionally, the accuracy is controllable over the omitted singular values.
The computation of statistics of the solution is straightforward. It is possible to compute multiple statistical scenarios with one single metamodel, since the information of probability distributions is used only to compute statistics and not to construct the metamodel.
\paragraph{Limitations}
If the random input parameters behave nonlinear on the solution, it is required to increase the number of sampling points (see line~\ref{algo:metamodel:simulate} of Algorithm~\ref{algo:accMetamodel}) to construct a metamodel with appropriate accuracy. This includes the execution of deterministic simulations on these sampling points to compute discrete solutions.

\section{Numerical Results}\label{sec:Results}
\subsection{Overview of the model problem}\label{sec:modelProblem}
\begin{figure}[htbp]
\centering
\begin{subfloat}[Geometry of the model problem]
{\psfrag{M1}[][]{\scriptsize{$\begin{array}{c}\mathrm{Material}\: 3 \\(\mathbf{D}_{3})\end{array}$}}
\psfrag{M2}[][]{\scriptsize{$\begin{array}{c}\mathrm{Material}\: 2 \\(\mathbf{D}_{2})\end{array}$}}
\psfrag{M3}[][]{\scriptsize{$\begin{array}{c}\mathrm{Material}\: 1 \\(\mathbf{D}_{1})\end{array}$}}
\psfrag{u = 0}[][]{\scriptsize{$u = 0$}}
\psfrag{u = 5}[][]{\scriptsize{$u = 1.0+\sigma_\mathrm{D} \xi_\mathrm{D}$}}
\psfrag{u.n = 0}[][]{\scriptsize{$\frac{\partial u}{\partial n} = 0$}}
\includegraphics[width=0.48\textwidth]{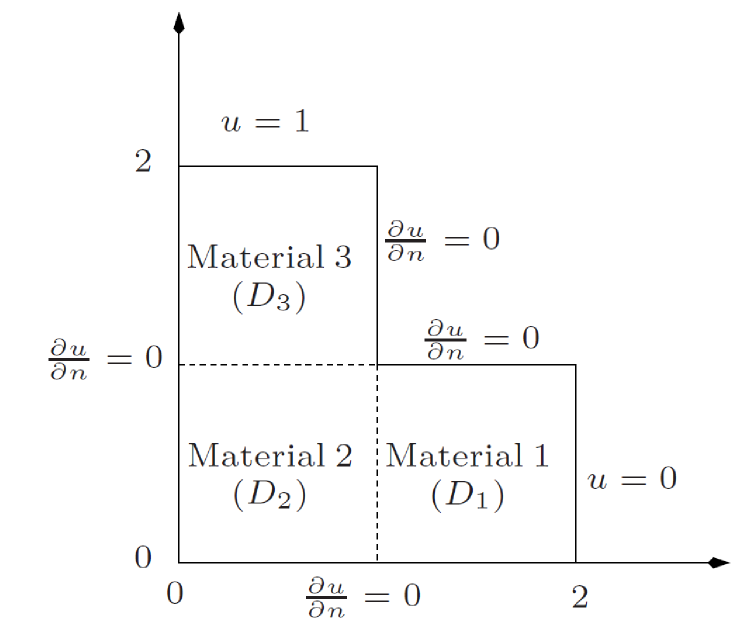}
}
\end{subfloat}
\begin{subfloat}[Deterministic solution]
{\includegraphics[width=0.45\textwidth,clip=true, trim=100mm 00mm 100mm 50mm]{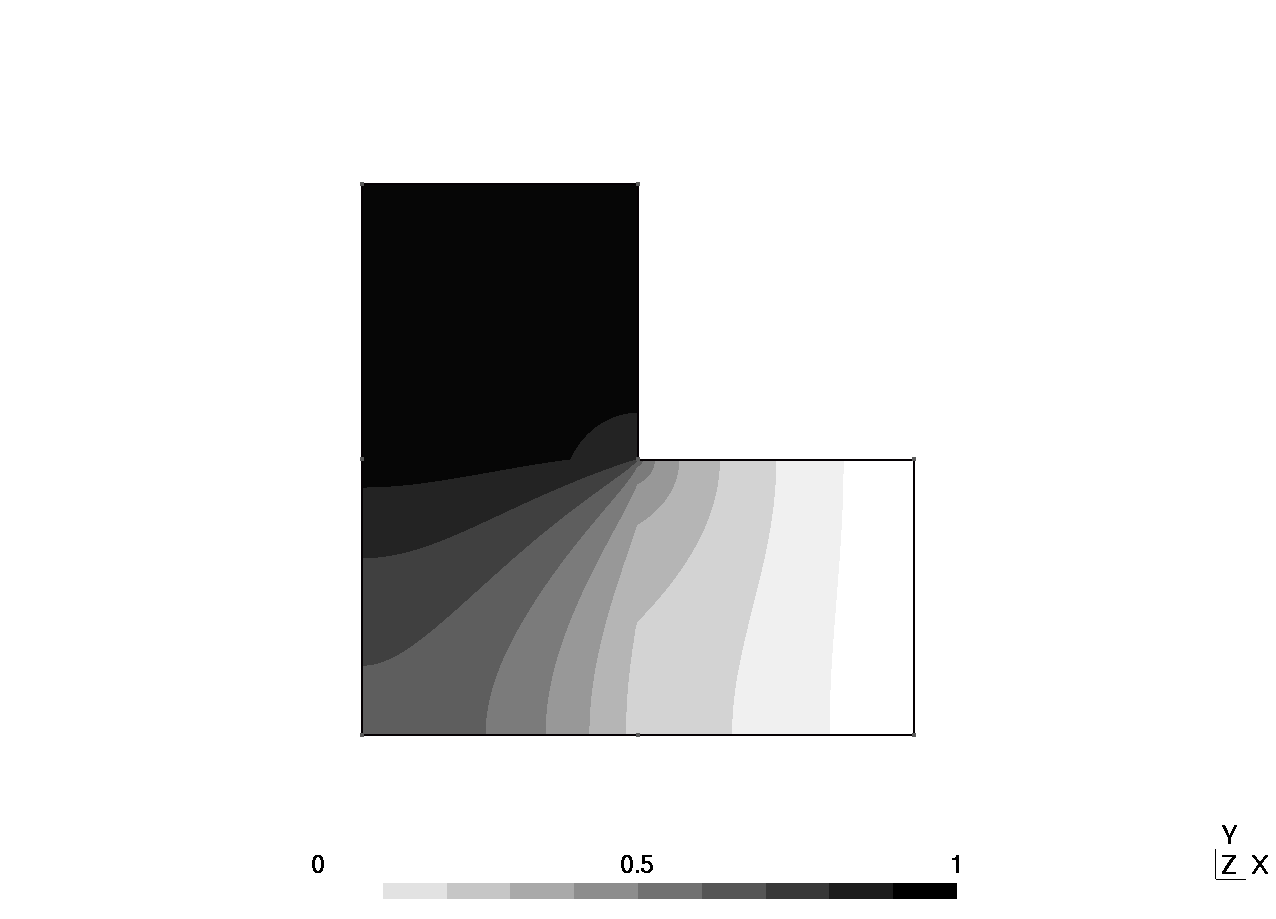}}
\end{subfloat}
\caption{Geometry and deterministic solution of model problem. In the problem shown here, all random fields $a_i$ are set to their mean values and $\gamma=1.0$. The problem is discretized on a mesh with 103479 nodes.}
\label{fig:Lshape}
\end{figure}

We consider a steady-state diffusion problem with a random diffusion coefficient, defined on an L-shaped domain~$\mathbf{D}$, as illustrated in Fig.~\ref{fig:Lshape},
\begin{align}\label{eq:diff1}
- \nabla \cdot (a(\boldsymbol{x},u(\boldsymbol{x},\boldsymbol{y}),\boldsymbol{y})\nabla u(\boldsymbol{x},\boldsymbol{y})) &= b(\boldsymbol{x}) \qquad \boldsymbol{x} \in \mathbf{D},\: \boldsymbol{y}\in \Gamma\:. 
\end{align}
Dirichlet boundary conditions are imposed on the upper and lower right boundaries, and zero Neumann conditions elsewhere. The right-hand side  $b(\boldsymbol{x})$ is equal to 1. We model $a(\boldsymbol{x},\boldsymbol{y})$ as a nonlinear piecewise random field, consisting of three parts, $a_1$, $a_2$ and $a_3$, defined respectively on the domains $\mathbf{D}_1$, $\mathbf{D}_2$ and $\mathbf{D}_3$:
\begin{equation}\label{eq:diffcoeff}
 a(\boldsymbol{x},u(\boldsymbol{x},\boldsymbol{y}),\boldsymbol{y}) =\left\{ \begin{array}{ccc}
                               a_1(\boldsymbol{x},\boldsymbol{y}) +\gamma u^2(\boldsymbol{x},\boldsymbol{y})& &\boldsymbol{x} \in \mathbf{D}_1\\
  a_2(\boldsymbol{x},\boldsymbol{y}) +\gamma u^2(\boldsymbol{x},\boldsymbol{y})& & \boldsymbol{x} \in \mathbf{D}_2\\
  a_3(\boldsymbol{x},\boldsymbol{y}) +\gamma u^2(\boldsymbol{x},\boldsymbol{y})& & \boldsymbol{x} \in \mathbf{D}_3
                                                            \end{array} \right. ,
\end{equation}
where $\gamma$ is a constant for varying the nonlinearity.
We assume that each component $a_i(\boldsymbol{x},\boldsymbol{y})$, $i = 1,2,3$, has the form of a truncated Karhunen-Lo\`eve (KL) expansion~\cite{Loeve1977} with respectively $L_1= 6$, $L_2 =7$ and $L_3 = 5$ random variables. We assume that all $L$ random variables are independent, with $L = L_1 + L_2 +L_3=18$. For the numerical experiments, we use an exponential covariance function for $a_i$, given by
\begin{equation}\label{eq:cov}
 C_{aa}(\boldsymbol{x},\boldsymbol{x}') = \sigma^2 \exp\left(-\frac{||\boldsymbol{x} - \boldsymbol{x}'||_1}{l_c}\right).
\end{equation}
The mean values of $a_i$ are set to $\bar{a}_1 = 30$, $\bar{a}_2 = 5$ and $\bar{a}_3 = 100$, respectively. 
Analytical expressions for the corresponding KL-functions are derived in \cite{Zhang2004}. We consider two configurations for creating the KL-expansions
of $a_i$ in \eqref{eq:diffcoeff}:
\begin{itemize}
\item $L$ uniformly distributed random variables on $[-\sqrt{3},\sqrt{3}]$: covariance function~\eqref{eq:cov} with correlation lengths $l_{c,1} = 1$, $l_{c,2} = 0.5$ and $l_{c,3} = 1.5$; and variances $\sigma_1^2 = 100$, $\sigma_2^2 = 2.25$, $\sigma_3^3 = 900$.
\item $L$ standard normally distributed random variables (mean $0$ and variance 1): covariance function~\eqref{eq:cov}
with $l_{c,1} = 1$, $l_{c,2} = 0.5$ and $l_{c,3} = 1.5$; and variances $\sigma_1^2 = 9$, $\sigma_2^2 = 0.25$, $\sigma_3^3 = 100$.
\end{itemize}
The problem is solved on a spatial finite element grid with $M=12154$ nodes and with $M=103479$ nodes. In each case, we set the constant $\gamma$ for varying the nonlinearity (\ref{eq:diffcoeff}) to $1.0$ and $100.0$. An example of a deterministic solution is shown in Figure~\ref{fig:Lshape}.

The nonlinear deterministic PDEs that are obtained by sampling the random variables are solved with a Newton-multigrid method, that is, algebraic multigrid is used to solve a sequence of linearized systems.

When applying the stochastic collocation method, a sparse grid of Clenshaw-Curtis or Gauss points is applied, denoted by `Scc' or `Sg', respectively. The level of the method will be indicated by the number following `Scc' and `Sg'.

Since the exact solution of the stochastic PDE is not known, we base the error computation on a reference solution $u_{\mathrm{ref}}$, which is obtained as a high-order stochastic collocation solution.
Then, we compare the accuracy of the different methods by computing the absolute difference $\mathrm{diff_{abs}}(j)$ or the relative difference $\mathrm{diff_{rel}}(j)$ locally in each grid point as
\begin{equation}\label{eq:errrel}
\mathrm{diff_{abs}}(j) = |\tilde{u}(j)-u_{\mathrm{ref}}(j)|, \qquad
\mathrm{diff_{rel}(j)} = \frac{|\tilde{u}(j)-u_{\mathrm{ref}}(j)|}{u_{\mathrm{ref}}(j)}, \quad j=1,\dots M, 
\end{equation}
where $\tilde{u}(j)$ denotes the interpolated result in grid point $j$ and $u_{\mathrm{ref}}(j)$ denotes the reference solution in grid point $j$.

\subsection{Parameter screening results}\label{sec:Results:screening}
We perform a design of experiments as described in Section \ref{sec:sensan:Doe} with $2L+1 = 37$ sampling points. Then, we apply the parameter screening measures of Section~\ref{sec:sensan:JH} to the model problem.

\begin{figure}[htbp]
\centering
\begin{subfloat}[$M=12154$, $\gamma = 1.0$]
{\includegraphics[width=0.48\textwidth]{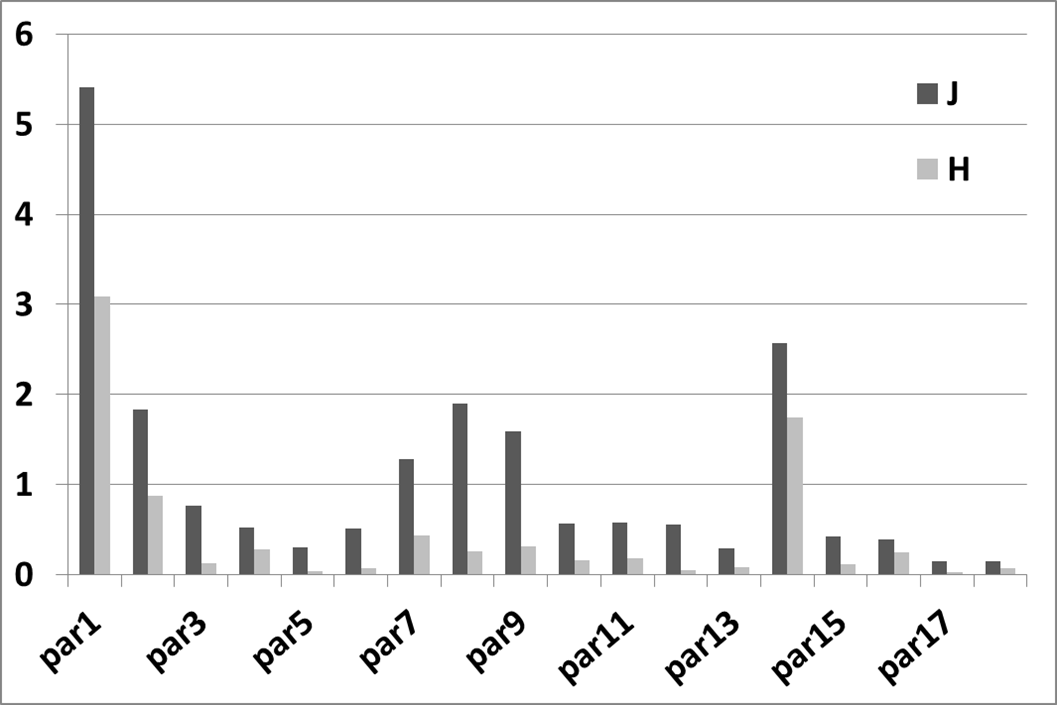}}
\end{subfloat}
\begin{subfloat}[$M=12154$, $\gamma = 100.0$]
{\includegraphics[width=0.48\textwidth]{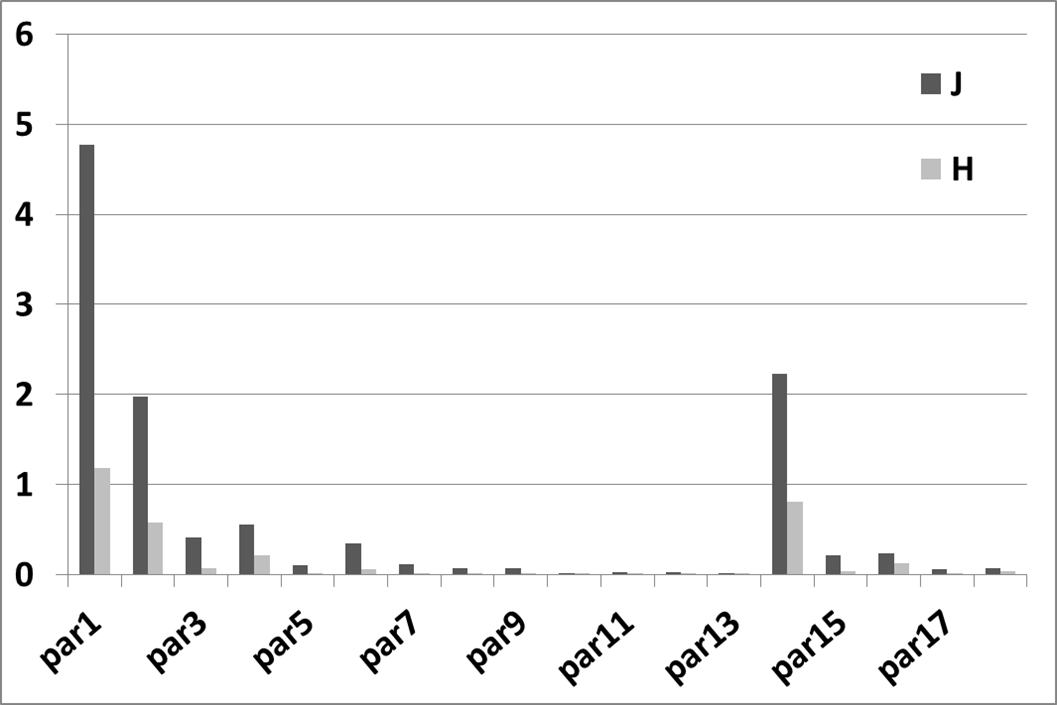}}
\end{subfloat} \\ 
\begin{subfloat}[$M=103479$, $\gamma = 1.0$]
{\includegraphics[width=0.48\textwidth]{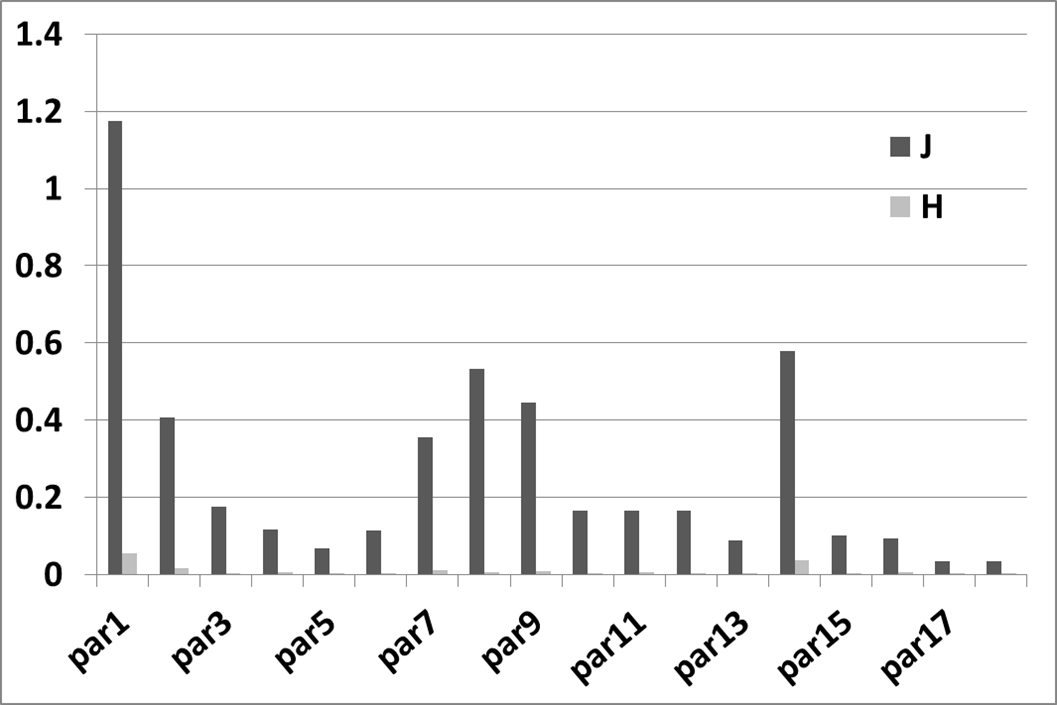}}
\end{subfloat}
\begin{subfloat}[$M=103479$, $\gamma = 100.0$]
{\includegraphics[width=0.48\textwidth]{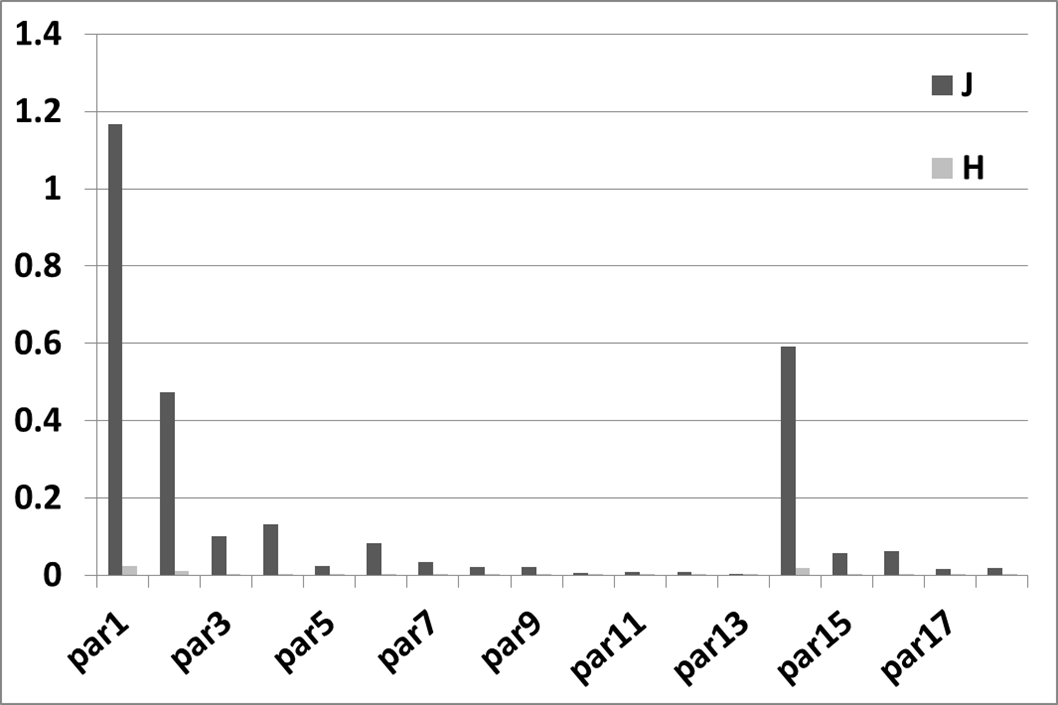}}
\end{subfloat}
\caption{Parameter screening results. The height of the bars gives the influence of the parameter on the solution, where the global approximations of $\mathbf{J}$ and $\mathbf{\mathrm{diag}H}$ are computed via \eqref{eq:sensitivityLinGlob} and \eqref{eq:sensitivityNonlinGlob} respectively.}
\label{fig:screening}
\end{figure}
The results of the parameter screening in the case of $M=12154$ nodes are shown in Figure~\ref{fig:screening}. The sensitivity analysis gives similar results for different resolutions of the grid. If we use only 12154 spatial degrees of freedom the indicator for nonlinearity \eqref{eq:sensitivityNonlinGlob} is higher than with a high resolution, but the ranking of the parameters remains the same. Differences in the ranking and the amount of influence of the parameters occur due to the variation of $\gamma$ in the model problem.
In detail, Figure~\ref{fig:screening} shows that in the case of $\gamma = 1.0$ the parameters $y_1,y_2,y_7,y_8,y_9$ and $y_{14}$ are the most influencing ones. The corresponding design-parameter space can be reduced to 6 instead of 18 parameters. In the case of $\gamma = 100.0$ the influence of the parameters at the first few ranks increases, so that the design-parameter space can be reduced to only three parameters, namely $y_1,y_2$ and $y_{14}$. The additional indicator for nonlinearity shows no great nonlinear influence of any parameter. In the cases where the measure \eqref{eq:criteriumGlob} is true, that is for $y_{10}$ and $y_{13}$ in the case of $\gamma=100.0$, we find that $J$ and $H$ are almost zero. Therefore these parameters do almost not influence the solution at all.

To confirm the linear behavior of all random variables, we compute the full Hessian matrix and the measure $D$~\eqref{eq:measureD} in case of 12154 nodes. Therefore, we have to perform $4L(L-1)/2 = 612$ additional simulations to compute the approximation of the second derivatives. We get
\begin{equation}
|\alpha_{\mathrm{max}}| = 1.05487, \quad D = 3.60796 > 0,
\end{equation}
with $\sigma = \sqrt{3}$. This proves the linear behavior of all random variables, at least up to second order approximations.  

\subsection{Comparison of the results: Collocation method and accelerated metamodel}
\begin{figure}[htbp]
 \centering
 \begin{subfloat}[Approximation error \eqref{eq:err} of the deviation $dX$ in dependence on $k$.]
{\includegraphics[width=0.45\textwidth]{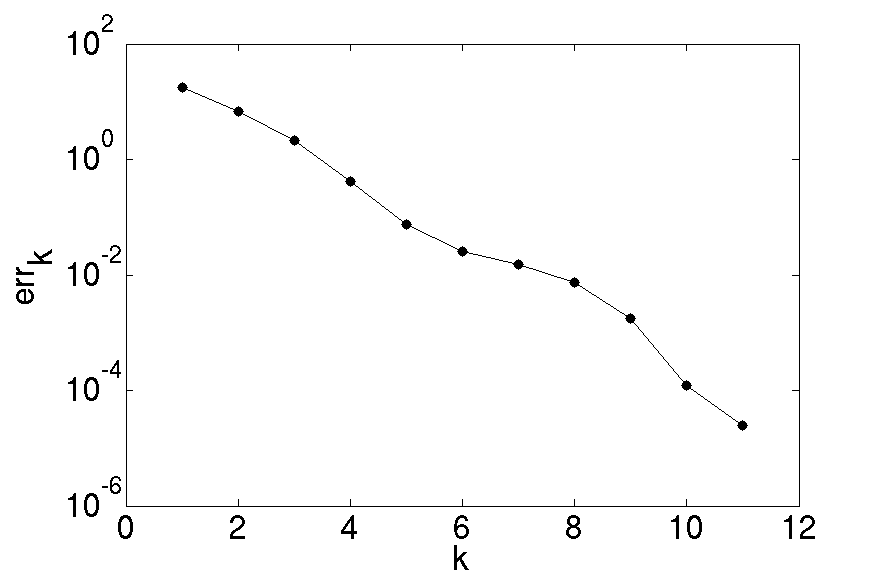}}
 \end{subfloat}
\begin{subfloat}[Corresponding singular values $\lambda_l$.]
{\includegraphics[width=0.45\textwidth]{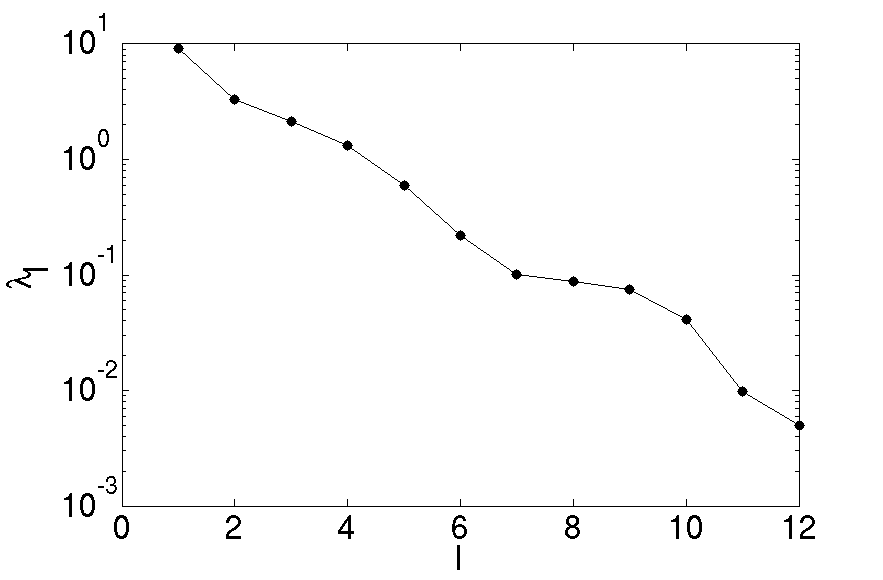}}
 \end{subfloat}\\
 \begin{subfloat}[Approximation error \eqref{eq:err} of the deviation $dX$ in dependence on $k$.]
{\includegraphics[width=0.45\textwidth]{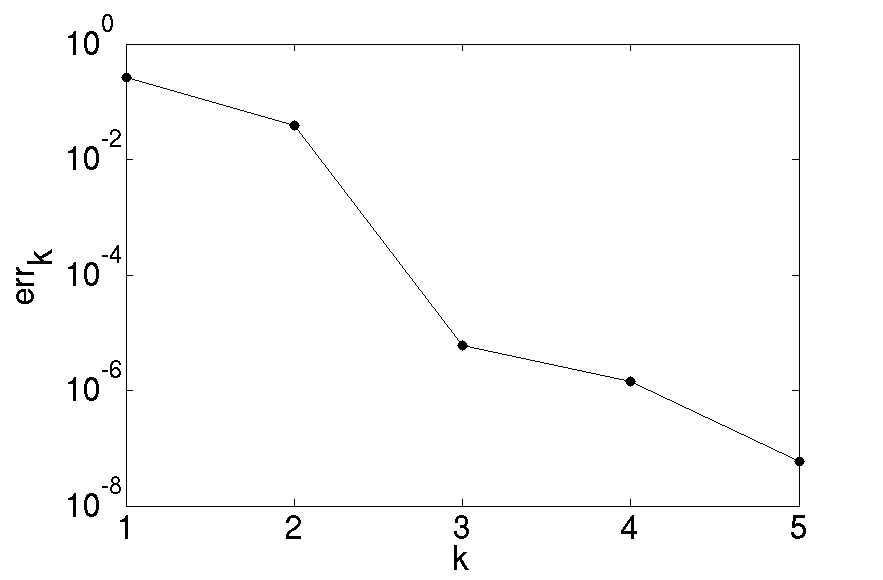}}
 \end{subfloat}
\begin{subfloat}[Corresponding singular values $\lambda_l$.]
{\includegraphics[width=0.45\textwidth]{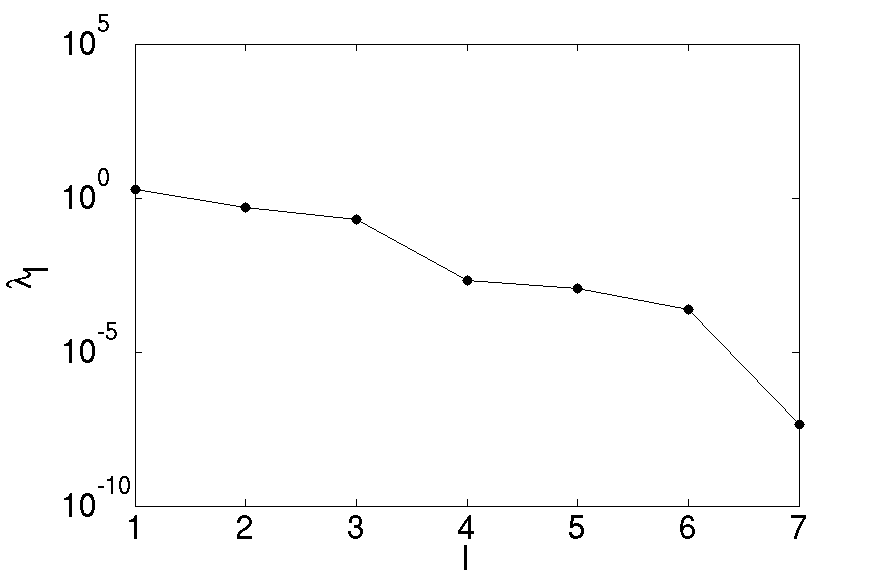}}
 \end{subfloat}
\caption{Precision of reconstructed matrix by \textit{fastSVD} as a function of $k$. (a)--(b)  18 uniformly distributed random variables, discretization with 12154 nodes, $\gamma=1.0$. (c)--(d) 18 Gaussian distributed random variables, discretization with 103479 nodes, $\gamma=100.0$.}
 \label{fig:parseval}
 \end{figure}
In order to show the efficiency of the accelerated metamodel, we compare the results to a reference solution, since the exact solution of the stochastic PDE is not known. We choose the collocation method as reference solution, since the collocation method is a usually applied state-of-the-art interpolation approach in the context of stochastic PDEs.

We compute the accelerated metamodel based on the multiquadric RBF, see Table~\ref{tab:rbf}, and $q$-quantiles according to Algorithm~\ref{algo:accMetamodel} for several configurations of the model problem~\eqref{eq:diff1}. In the following, two cases are presented in detail.
First, we analyse the case of 18 uniformly distributed random variables on a mesh with 12154 nodes and nonlinearity coefficient $\gamma=1.0$. Here, the design-parameter space is reduced to only 6 random variables as a result of the previous parameter screening, see Section~\ref{sec:Results:screening}. Therefore, the sample points for the construction of the metamodel also reduce to $2L+1 = 13$. We set up the accelerated metamodel for this configuration with 6 random variables.

As a second example, we analyse the case of 18 Gaussian distributed random variables on a mesh with 103479 nodes and nonlinearity coefficient $\gamma=100.0$. Here, the design-parameter space is reduced to 3 random variables as a result of the previous parameter screening and $7$ sample points are used for the construction of the accelerated metamodel. In the following, we will always use the metamodels with the reduced design-parameter space, unless stated otherwise.

\subsubsection{RBF metamodel}
We compute the \textit{fastSVD} for each model problem as described in Section~\ref{sec:svd}. The precision of the reconstructed matrix as a function of $k$, together with the corresponding singular values $\lambda_l$, is illustrated in Fig.~\ref{fig:parseval} exemplary for a discretization with 103479 nodes and $\gamma=100.0$. The figures show that the singular values are rapidly decreasing with $k$. Therefore, the error is rapidly decreasing with $k$ as well. Only 7 to 10 singular values in the case of uniformly distributed random variables and 3 singular values in the case of Gaussian distributed random variables are sufficient to get an accurate solution.

Next, the metamodel weights are computed according to~\eqref{eq:weights}.
To evaluate this metamodel, we perform a new DoE with $L_{add} = 2000$ sampling points $\tau_j, j=1,\dots,L_{add}$. The sampling points are randomly chosen and uniformly distributed in the design space.
The deterministic solutions of the sPDE according to these sampling points are computed by interpolation with the accelerated metamodel, using formula \eqref{eq:mixing}.
\begin{figure}[htbp]
 \centering
 \begin{subfloat}[Deterministic solution for sampling point $\tau_1$.]
{\includegraphics[width=0.45\textwidth,clip=true, trim=100mm 00mm 80mm 50mm]{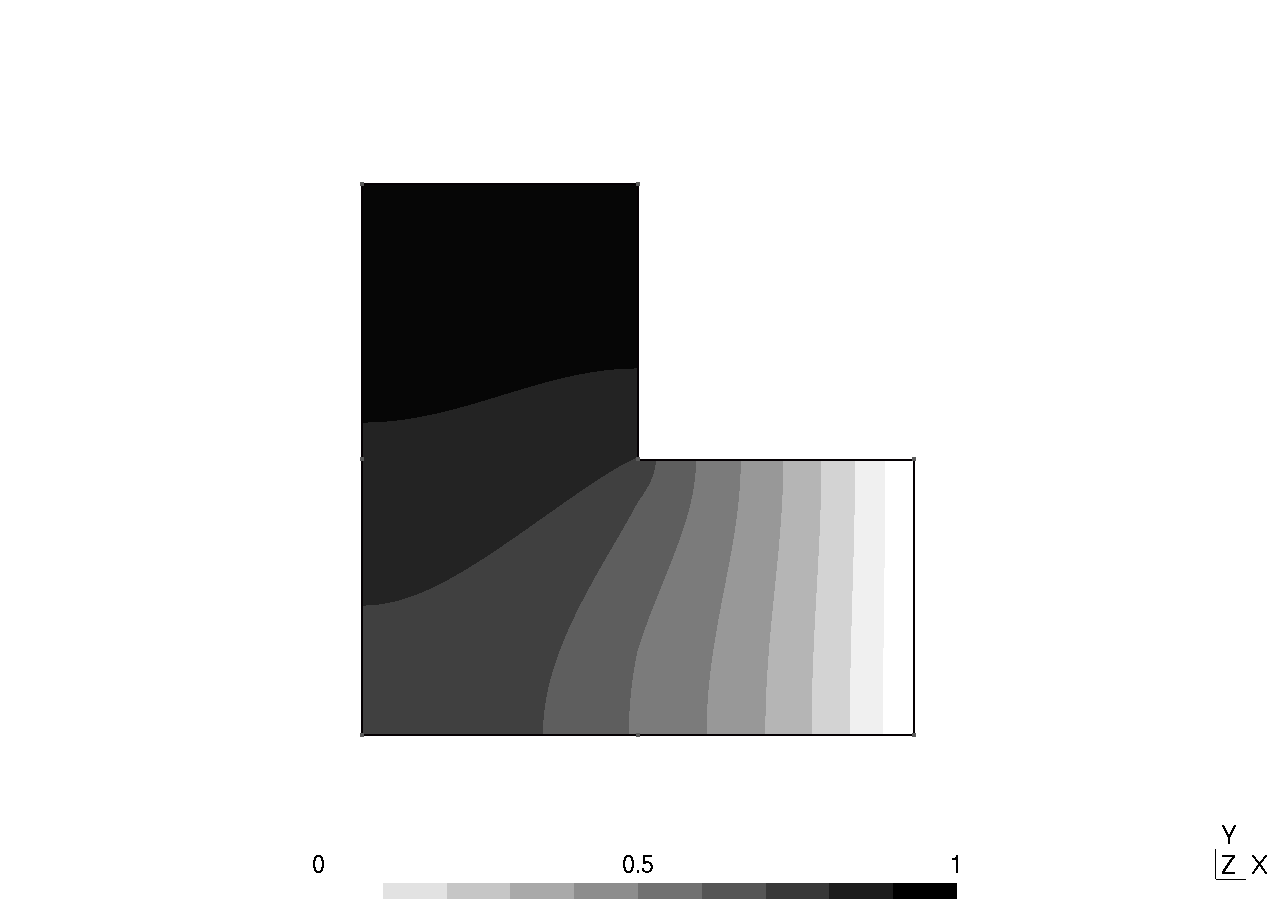}}
 \end{subfloat}
\begin{subfloat}[Absolute difference between the interpolated result in $\tau_1$ and the deterministic solution.]
{\includegraphics[width=0.45\textwidth,clip=true, trim=100mm 00mm 80mm 50mm]{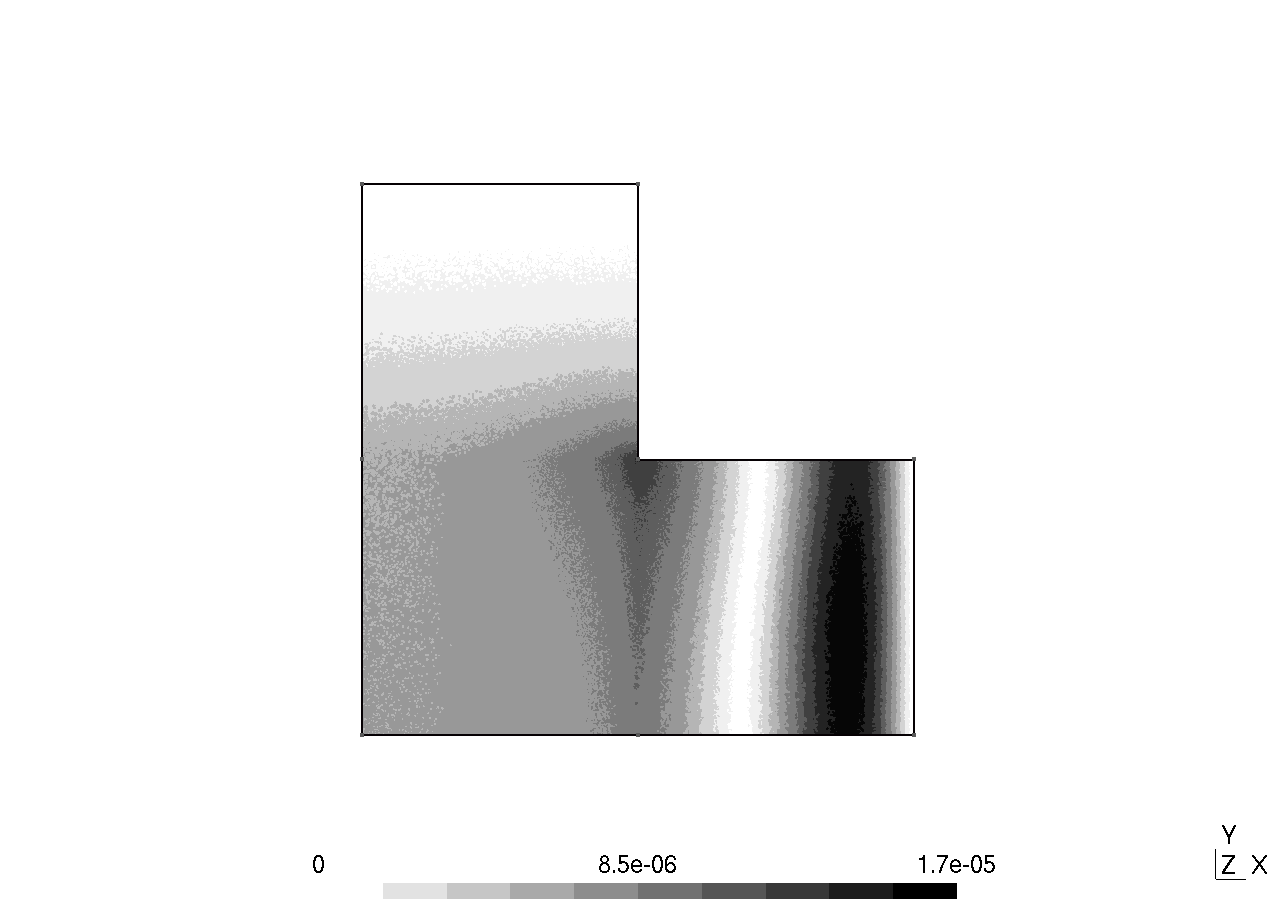}}
 \end{subfloat}
\caption{Deterministic solution and absolute difference between the interpolation result (reduced metamodel with 3 random variables) and this deterministic solution in the case of 18 Gaussian distributed random variables, using $\gamma=100$ and a spatial discretization of 103479 nodes.}
 \label{fig:absdiff-tau}
 \end{figure}
As an example, the deterministic result for one single sample point $\tau_1 =[  
\xi_1 = 0.2, \xi_2 = 0.4, \xi_{14} = 0.8, \xi_i=0.5 \quad \forall i \backslash \{1,2,14\} ] $ is shown in Figure~\ref{fig:absdiff-tau} on the left, on the right side the absolute difference between the interpolation result in this sample point and this deterministic solution is shown. Note that for the interpolation result the metamodel with the reduced parameter space has been used. Figure~\ref{fig:absdiff-tau} shows that the interpolation with the accelerated metamodel gives accurate results. With only the reduced parameter space the maximal absolute error is $1.72907e-05$. The maximal relative error is $0.000122$, that is $0.012$ percent of the solution value, and this appears only in the region, where the deterministic solution equals almost zero.

The whole set of sample points and corresponding interpolation results is then used to compute several quantiles, for example the median, the $0.68$-quantile and the $0.9$-quantile.

\subsubsection{Collocation method versus accelerated metamodel}
We compare the DoE and the interpolation methods with respect to the number of deterministic PDEs to be solved, the computational time and accuracy.
The collocation methods and the accelerated metamodel require a different number of sampling points, see Table~\ref{tab:NSamplePoints}, and therewith a different number of deterministic PDEs to solve.
\begin{table}[htbp]
  \centering
  \caption{Number of required sampling points as a function of the number of random variables}
    \begin{tabular}{|r|rrrr|r|}
    \addlinespace
    \toprule
    \textbf{random variables} & \textbf{Scc} & \textbf{Scc2} & \textbf{Sg} & \textbf{Sg2} & \textbf{accMetamodel} \\
    \midrule
    18    & 37    & 685   & 37    & 757   & 37 \\
    6     & 13    &  85     & 13    & 109   & 13 \\
    3     & 7     & 25    & 7     & 37    & 7 \\
    \bottomrule
    \end{tabular}%
  \label{tab:NSamplePoints}%
\end{table}%

Details of the different components of the computational time are shown in the Tables~\ref{tab:time-coll-103479} and~\ref{tab:time-metamodel-103479}.

\begin{table}[htbp]
  \centering
  \caption{Detailed overview of runtime (in sec) on a standard 3GHz Linux PC to compute $0.68$-quantile with the stochastic collocation method for the model problem with 18 Gaussian random variables, $\gamma=100.0$ and a spatial discretization with 103479 nodes.}
    \begin{tabular}{|r|rr|}
    \addlinespace
    \toprule
    \textbf{Time for} & \textbf{Sg} & \textbf{Sg2} \\
    \midrule
    construction of collocation points & 0.000000E+00 & 0.000000E+00 \\
    solution of deterministic PDEs (simulation) & 2.171000E+03 & 4.273300E+04 \\
    construction of samples of the collocation solution & 1.919000E+01 & 4.162900E+02 \\
    0.68-quantile & 1.022000E+01 & 1.044000E+01 \\
    \textbf{overall runtime} & \textbf{2.200410E+03} & \textbf{4.315973E+04} \\
    \bottomrule
    \end{tabular}%
  \label{tab:time-coll-103479}%
\end{table}%
\begin{table}[htbp]
  \centering
  \caption{Detailed overview of runtime (in sec) on a standard 3GHz Linux PC to compute $0.68$-quantile with the accelerated metamodel method for the model problem with 18 Gaussian random variables, $\gamma=100.0$ and a spatial discretization with 103479 nodes.}
    \begin{tabular}{|r|r|}
    \addlinespace
    \toprule
    \textbf{Time for} & \textbf{acc\_metamodel\_18} \\
    \midrule
    construction of DoE for 18 variables & 0.000000E+00 \\
    solution of deterministic PDEs (simulation) & 2.171000E+03 \\
    parameter screening & 3.504000E+00 \\
    construction of DoE for 3 variables & 0.000000E+00 \\
    fastSVD & 1.370000E+00 \\
    evaluation of the accelerated metamodel and 0.68-quantile & 1.121000E+01 \\
    \textbf{overall runtime} & \textbf{2.187084E+03} \\
    \bottomrule
    \end{tabular}%
  \label{tab:time-metamodel-103479}%
\end{table}%
Beside the different computational time for the simulations, the different times for constructing the sampling points for the statistics are standing out. The collocation method uses $2000$ sampling points according to the probability distribution of the random variables to reconstruct the stochastic solution~\eqref{eq:lagrange}. We use the same number of sampling points $L_{add} = 2000$ also for the evaluation of the accelerated metamodel.
The tables demonstrate that the parameter screening is very fast, although it computes measures in each grid point. The time for the evaluation of the accelerated metamodel and the 0.68-quantile contains the time for computing the solution of the sPDE with interpolation in the current sampling point using~\eqref{eq:mixing}, added up for all sampling points and the update of the quantile estimator with the $P^2$ algorithm. Note that we use the same algorithm for the computation of quantiles, namely the $P^2$algorithm to get a meaningful comparison. Therefore, the higher runtime to compute the $0.68$-quantile with the collocation method is only due to the reconstruction of the stochastic solution.
Altogether, the runtime of the accelerated metamodel approach is all about the same than the collocation method with Gauss points, level one (Sg), and a lot faster than the collocation method with Gauss points, level two (Sg2).
Note that this comparison only illustrates where the computational effort of each method lies on. Each method may be accelerated by parallelizing components of the method.
\begin{figure}[htbp]
 \centering
 \begin{subfloat}[$0.68$-quantile, Metamodel vs Sg.]
{\includegraphics[width=0.4\textwidth,clip=true, trim=100mm 00mm 100mm 50mm]{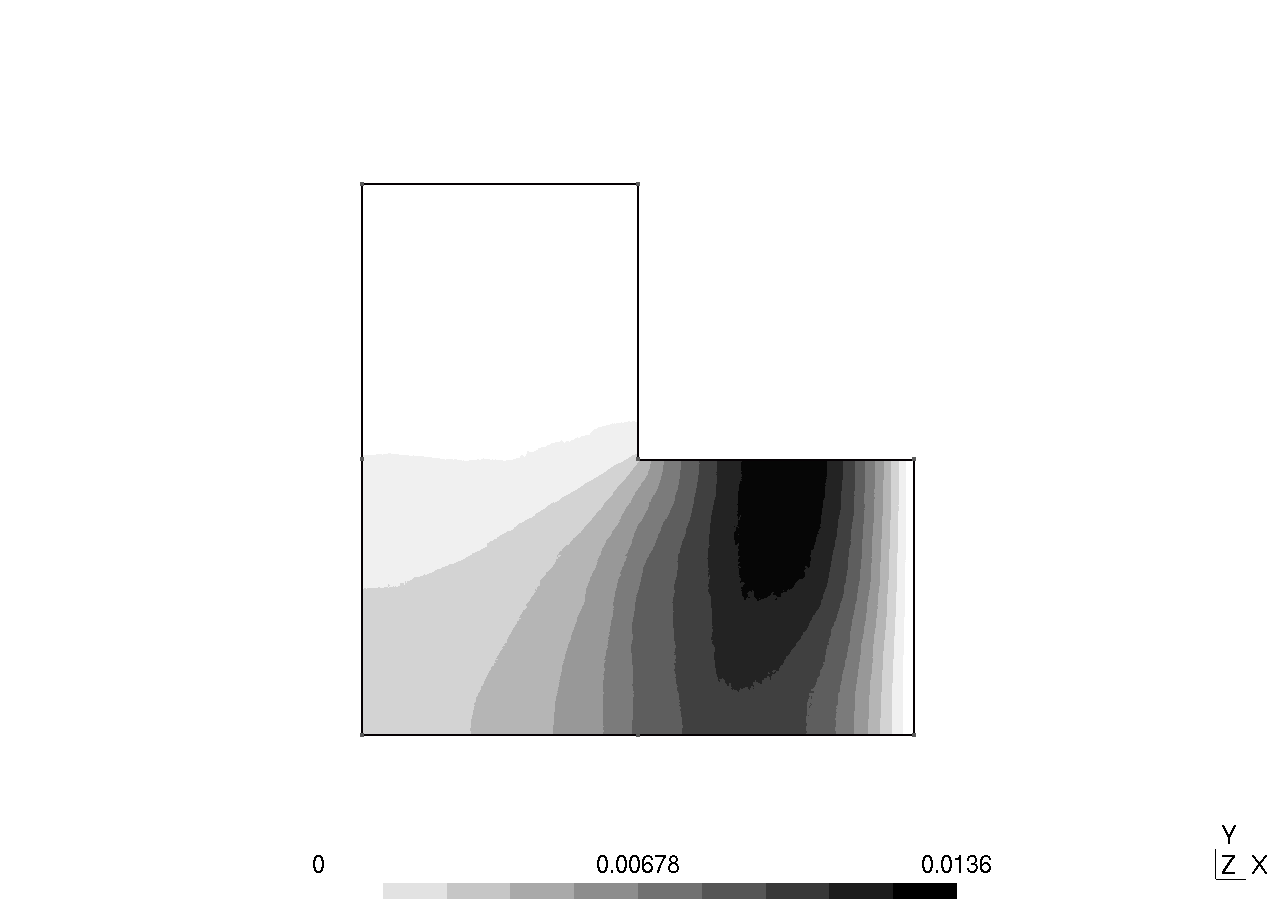}}
 \end{subfloat}
 \begin{subfloat}[$0.68$-quantile, Metamodel vs Sg2.]
{\includegraphics[width=0.4\textwidth,clip=true, trim=100mm 00mm 100mm 50mm]{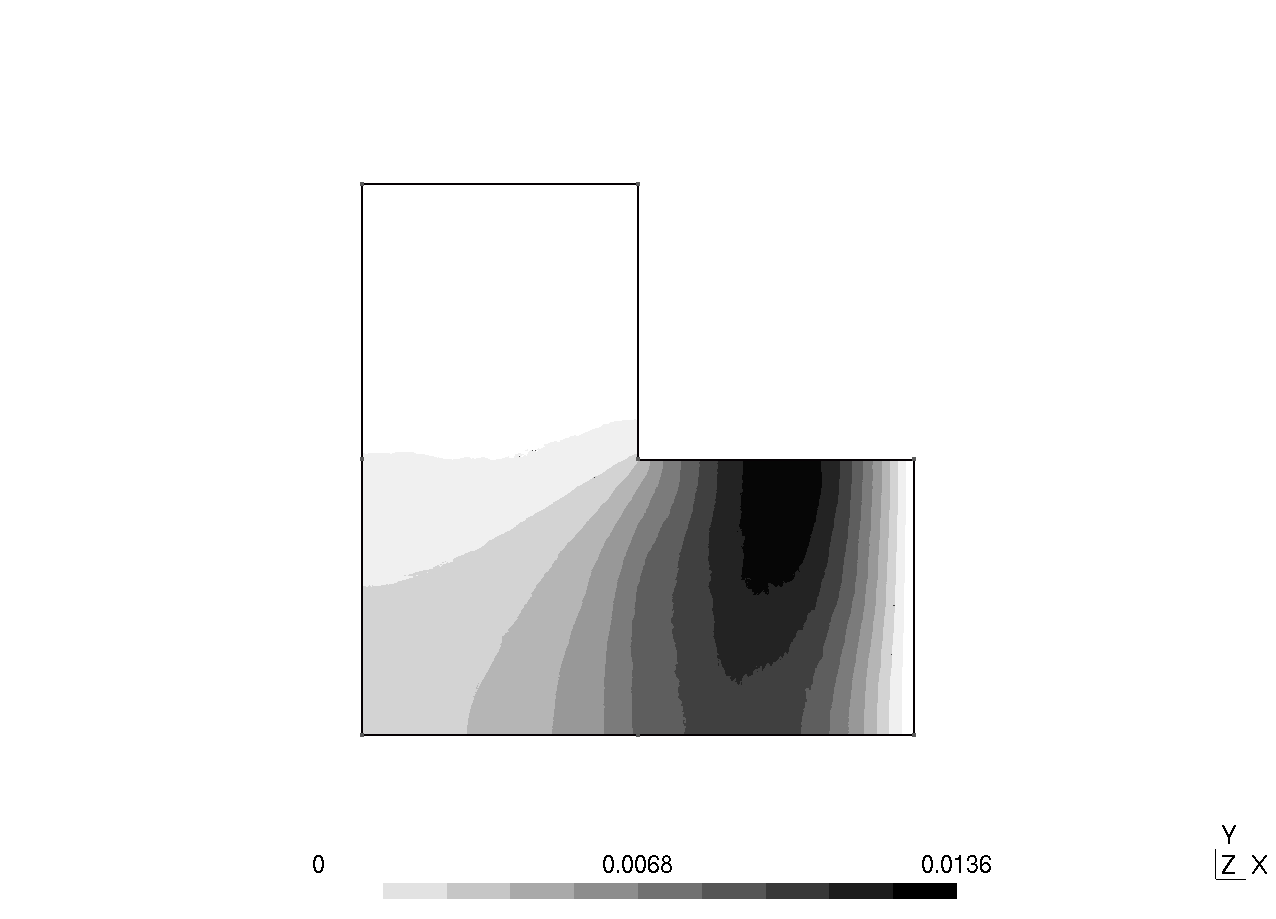}}
 \end{subfloat}
\caption{Absolute difference between the $0.68$-quantile obtained with the collocation and accelerated metamodel method. The model problem with 18 Gaussian distributed random variables, discretized with 103479 nodes, $\gamma=100.0$ and a reduced metamodel with 3 random variables are used.}
 \label{fig:quantile068-abs}
 \end{figure}
\begin{figure}[htbp]
 \centering
 \begin{subfloat}[$0.68$-quantile, Metamodel vs Sg.]
{\includegraphics[width=0.4\textwidth,clip=true, trim=100mm 00mm 100mm 50mm]{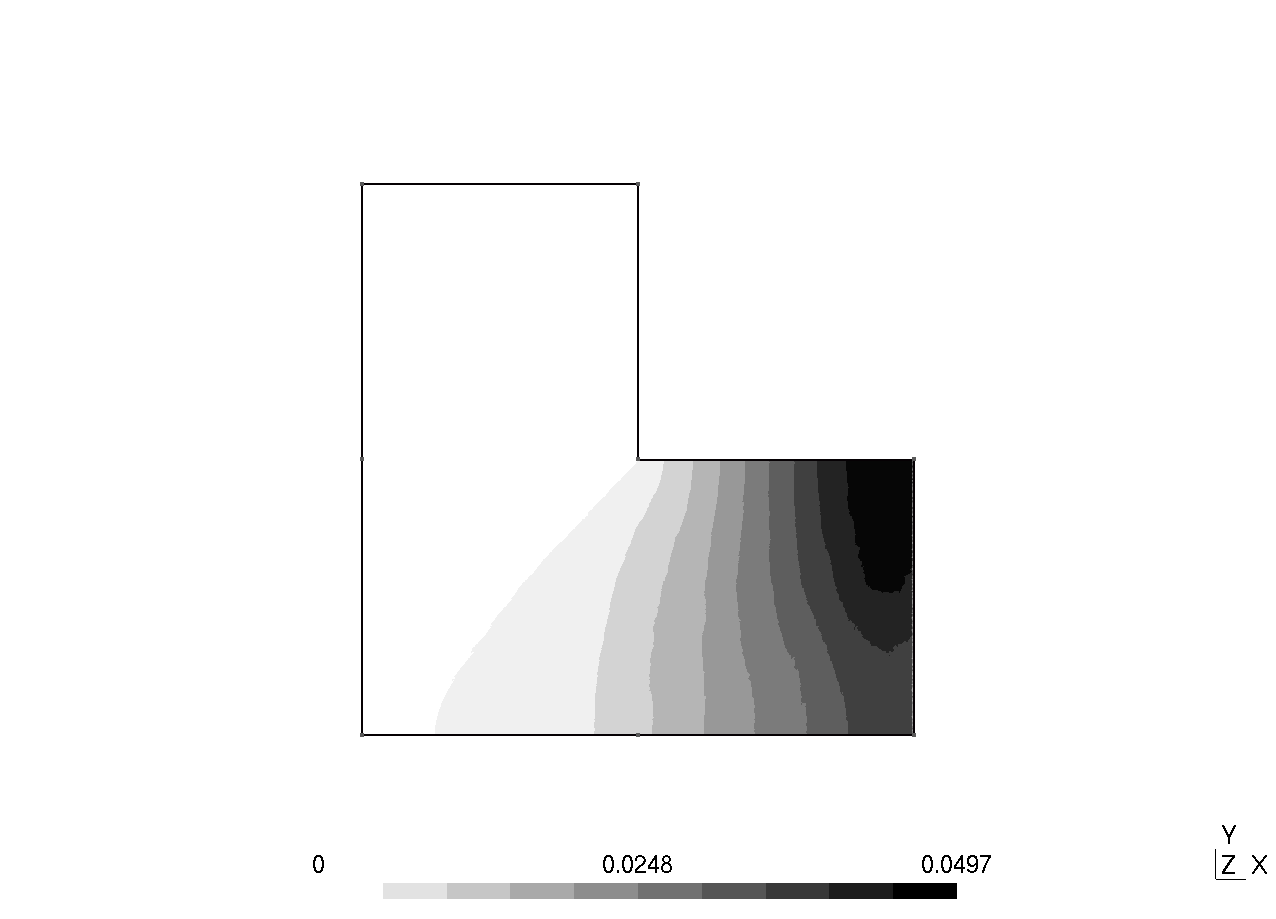}}
 \end{subfloat}
 \begin{subfloat}[$0.68$-quantile, Metamodel vs Sg2.]
{\includegraphics[width=0.4\textwidth,clip=true, trim=100mm 00mm 100mm 50mm]{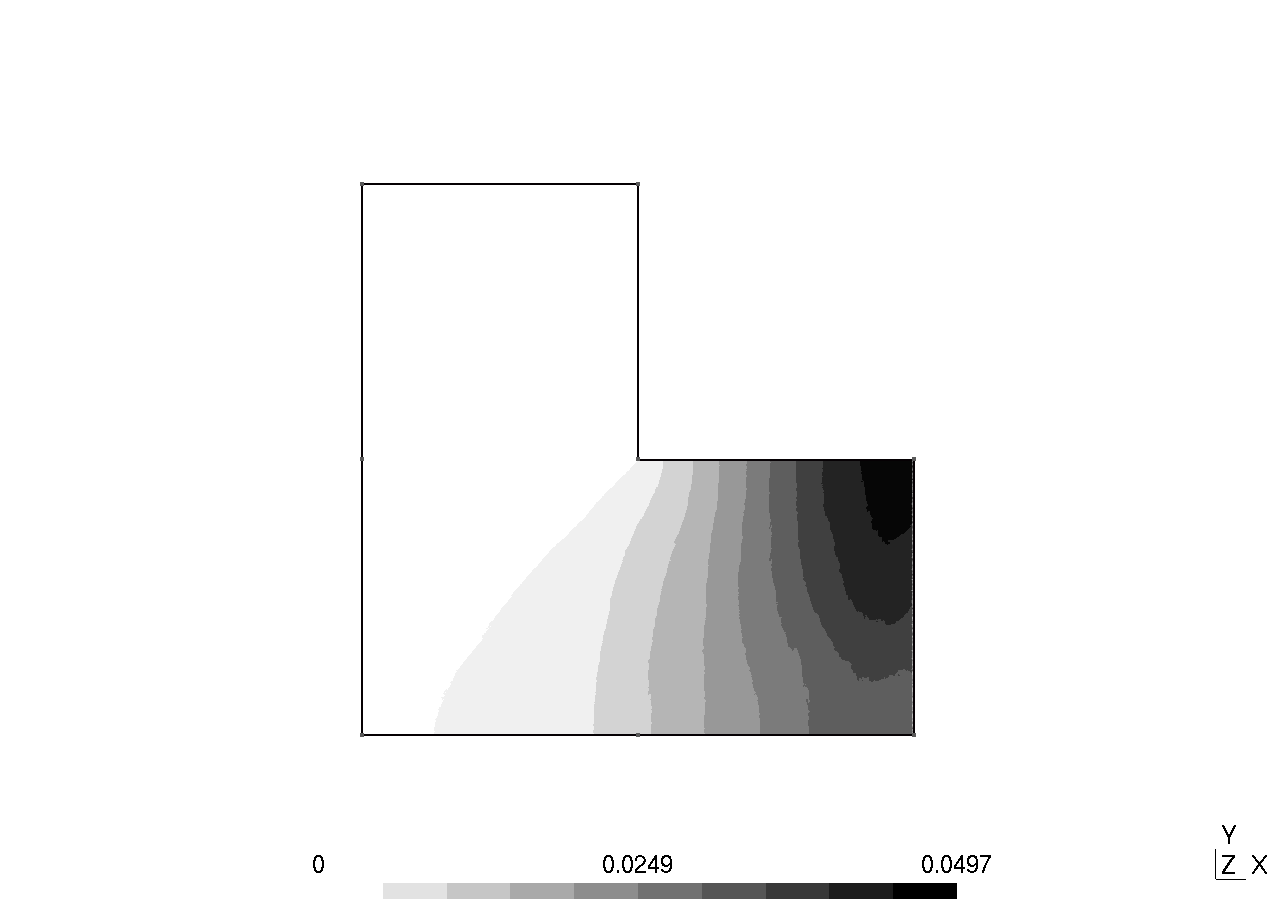}}
 \end{subfloat}
\caption{Relative difference between $0.68$-quantile obtained with the  collocation and accelerated metamodel method. The model problem with 18 Gaussian distributed random variables, discretized with 103479 nodes, $\gamma=100.0$ and a reduced metamodel with 3 random variables are used.}
 \label{fig:quantile068-rel}
 \end{figure}

We compare the metamodel solution against reference solutions $u_{ref} = Sg$ and $u_{ref} = Sg2$. However, note that the computed quantile with $Sg$ or $Sg2$ is the exact solution of the constructed polynomial solution approximation, but not the exact solution of the sPDE.
Figures~\ref{fig:quantile068-abs} and \ref{fig:quantile068-rel} show the absolute and relative differences~\eqref{eq:errrel} between the $0.68$-quantile computed with the collocation and with the accelerated metamodel method. Some differences occur, especially in the domain where the solution of the sPDE is almost zero. The differences are smaller when compared against $Sg2$. The errors are of course influenced by the number of sampling points. For comparison each method uses 2000 sampling points, but an evaluation of the accelerated metamodel with more samples may even ameliorate the result.

In summary, the accelerated metamodel performs as fast as the collocation Sg method for the model problem considered. But, its accuracy lies somewhere between Sg and Sg2.
Therefore, the numerical results illustrate that the accelerated metamodel approach is an efficient and accurate alternative to the collocation method to compute statistics of a sPDE containing independent random variables. It is especially advantageous when the probability distributions of the random variables are not known a priori.

\section{Discussion and future work}\label{sec:concl}
The paper deals with the computation of statistics of a sPDE by means of interpolation methods. We present an alternative approach to the collocation method, namely the accelerated RBF metamodel. The different methods were demonstrated on a nontrivial model problem. Numerical examples show the high potential of this new interpolation approach, which is very fast and needs only a small set of deterministic solutions of the sPDE. The great advantage of the presented method is its flexibility, it makes only few assumptions and does not require the exact probability distribution of the random input variables for construction. It is possible to compute multiple statistical scenarios with one single metamodel. Therefore, the approach will be applicable in most industrially relevant applications. 

Secondly, this paper deals with the curse of dimensionality. The presented parameter screening approach can be used in order to reduce the model complexity. This generally accelerates the computation of the solution of the sPDE. Here, the parameter screening approach is used with the accelerated metamodel approach. But, it is also possible to combine the parameter screening with the collocation method to get an accurate solution in less computational time. This will result in the construction of an adaptive sparse collocation method.

Future work will be on extending the parameter screening and accelerated metamodel approach to more nonlinear cases. This will lead to a better approximation of the deterministic solution and, thus, the $q$-quantiles.
Also, there is some work to do in order to find the optimal number of sampling points required for the computation of a $q$-quantile.
 
\bibliographystyle{abbrv}
\bibliography{references}
\end{document}